\documentclass[a4paper,10pt,epsf]{amsart}

\usepackage{amssymb,amsmath,amsthm}
\usepackage{latexsym}
\usepackage[english]{babel} 
\usepackage[latin1]{inputenc}  
\usepackage[T1]{fontenc}   
\usepackage[all]{xy}
\usepackage{multicol}
\usepackage{amsfonts}
\usepackage{graphics}
\usepackage{epsfig}
\usepackage{graphicx}

\newcommand{\color}[6]{}

\def\RR{\mathbb{R}}
\def\CC{\mathbb{C}}

\def\NN{\mathbb{N}}
\def\QQ{\mathbb{Q}}
\def\PP{\mathbb{P}}
\def\SS{\mathbb{S}}

\def\gG{\mathfrak{g}}

\def\a{\alpha}

\def\g{\gamma}

\def\onto{\longmapsto}

\newtheorem{prop}{Proposition}

\newtheorem{lemm}{Lemma       }
\newtheorem{def }{Definition  }[section]
\newtheorem{theo}{Theorem     }
\newtheorem{ex}{Example     }[section]
\newtheorem{rem}{Remark    }[section]
\newtheorem{coro}{Corollary     }[section]

\begin{document}

\title{Singularities of Duals of Grassmannians}

\author{Fr\'ed\'eric Holweck}

\maketitle

\begin{abstract}
 Let $X\subset \PP^{N}$ be a smooth irreducible nondegenerate projective variety
and let $X^*\subset \PP^{N}$ denote its dual variety. %It is well known that $\sigma_2(X)^*$, the dual 
%of the 2-secant variety of $X$, is a component of the singular locus of $X^*$. 
The locus of bitangent hyperplanes, i.e. hyperplanes tangent to at least two points of $X$, is a component of the sigular locus of $X^*$.
In this paper we provide a sufficient condition for this
 component to be of maximal dimension and show how it can be used to determine which 
dual varieties of Grassmannians are normal. That last part may be compared to what has been done for hyperdeterminants
by J. Weyman and A. Zelevinski (1996) in \cite{WZ}.
\end{abstract} 

\section{Introduction}
Let $X=\PP^n\times \PP^n \subset \PP^{(n+1)^2-1}$ be the Segre
embedding of the product of two projective spaces of dimension
$n$. The variety $X$ corresponds to the projectivization of the
variety of rank one matrices embedded in the projectivization of the
space of $(n+1)\times(n+1)$ matrices. It is well known its  dual variety (the
variety of tangent hyperplanes, see below for the definition), denoted
by $X^*$, can be identified with the variety of rank at most $n$ matrices. Up to
multiplication by a nonzero scalar, the equation defining $X^*\subset
\PP^{(n+1)\times(n+1)-1 *}$ is the 
determinant. That point leads  to  a higher dimensional
generalization of the determinant, called {\em hyperdeterminant}, which was first introduced by Cayley (1840) and 
rediscovered by Gelfand, Kapranov and Zelevinsky (1992). 
In \cite{GKZ1,GKZ2} the authors define the
 hyperdeterminant of format $(k_1+1)\times\dots\times (k_r+1)$ by the
equation (up to scale) of the dual variety of $X=\PP^{k_1}\times
\dots\times \PP^{k_r} \subset
\PP^{(k_1+1)\times\dots\times(k_r+1)-1}$. When the dual variety $X^*$
is not a hypersurface, the corresponding hyperdeterminant is defined to be zero.

Let  $X\subset \PP(V)$ be a projective variety and let $\tilde{T}_x
  X$ denote the embedded tangent space of $X$ at a point $x\in
  \text{Sm}(X)$ (smooth points of $X$). Define the {\em dual variety}  $X^*$ by

$$X^*=\overline{\{H \in \PP(V^*) | \ \exists \ x \in \text{Sm}(X) \text{ such that} \  \tilde{T}_x X\subset H\}}\subset\PP(V^*).$$

The biduality theorem $(X^*)^*=X$ (true in characteristic zero)
implies that the original variety can be reconstructed from its dual variety. Thus geometric
invariants of $X^*$ reflect in geometric properties of $X$. The dimension, degree and singularities of $X^*$ carry
 meaningful information about the hyperplane sections of $X$ (see
\cite{Z5}).
These invariants have been studied for hyperdeterminants. In
\cite{GKZ1} a condition is given to decide whether or not the
hyperdeterminant of a given format is nonzero (i.e. the dual of the Segre
embedding is actually a hypersurface), moreover in the same paper the authors give a
combinatorial formula to compute the degree of a given
hyperdeterminant. They also conjectured that there is only one
hyperdeterminant whose corresponding hypersurface is regular in
codimension one, i.e. $\text{codim}_{X^*} \text{Sing}(X^*)\geq 2$, and this hyperdeterminant is of format $(2,2,2)$. In
other words $\PP^1\times \PP^1\times \PP^1\subset \PP^7$ is the only
Segre  product of at least three projective spaces whose dual variety is a normal
hypersurface. That conjecture was proved
by Weyman and Zelevinsky in \cite{WZ}. 

Let $G(k,n)\subset \PP^{\binom{n}{k}-1}$ denote the Grassmannian of $k$-planes in $V=\CC^n$, $k\leq n-k$, embedded through the Pl\"ucker map. Its dual variety is a hypersurface except if $k=2$ and $n$ is odd \cite{KM}. The degree of $G(k,n)^*$ has been studied in \cite{Las}. 
However the study of $\text{Sing}(G(k,n)^*)$ has not been carried out so far. In this article we answer the question of the normality of the duals of Grassmannian varieties. The case of the Grassmannian of
$2$-planes is known and similar to the Segre product of two projective spaces. The variety $G(2,n)\subset \PP(\Lambda^2 \CC^n)$ corresponds to the projectivization of the rank $2$ skew-symmetric matrices and its dual is identified with 
degenerate skew-symmetric matrices. Like for the determinant, the singular locus of the degenerate skew-symmetric matrices is regular in codimension $1$ and arithmetically Cohen Macaulay (\cite{Ke}).
This proves that $G(2,n)^*$ is normal. 

For $k\geq 3$ the dual variety $G(k,n)^*$ is a hypersuface. Thus $G(k,n)^*$ will be normal if and only if $G(k,n)^*$ is regular in codimension one.
This will be the main result of this article:

\begin{theo}\label{maintheo}
 Let $X=G(k,n)\subset \PP^{\binom{n}{k}-1}$, with $k\geq 3$. The dual variety $G(k,n)^*$ is normal if and only if $X$ is one of the following : 
\[G(3,6)\subset \PP^{19}, \ G(3,7)\subset \PP^{34}, \ G(3,8)\subset \PP^{55}\]
\end{theo}

\begin{rem}
 Like for hyperdeterminants the general pattern is the following: 
the variety $X^*$ has a singular locus of codimension one and the only exceptions
 come from group actions with finite numbers of orbits.
\end{rem}

The proof is based on the calculation of the dimension of $\sigma_2(X)^*$, the dual of the secant variety of $X$,  which is always a component of $\text{Sing}(X^*)$. It turns out that this component appears in the 
decomposition of the singular locus of hyperdeterminants by \cite{WZ}.
In their paper it corresponds to the general double point locus or node locus denoted by $\nabla_{node}(\emptyset)$ 
(i.e. the set of hyperplane having more than one point of tangency on $X$). An other component of interest is the cusp locus 
(i.e. set of hyperplanes defining degenerate quadrics). 
The geometrical meaning of $\nabla_{node}(\emptyset)$ is not emphasized in \cite{WZ} 
when they calculate  the dimension of this component.
Here in the contrary we mainly use geometric arguments to calculate the dimension of $\sigma_2(X)^*$ 
in the general case.
Let $\hat{T}_x^{(2)}X$ be the (cone over the) second osculating space, i.e. the linear span of 
second osculating spaces of smooth 
curves $x(t)\subset X$ with $x(0)=x$. In $\S\ref{comp}$ we prove:
\begin{prop}~\label{crit}
Let $X\subset \PP(V)$ be a smooth projective variety of dimension $n$.
Assume $X^*$ is a hypersurface. Suppose for a general pair of point $(x,y)\in X\times
X$ we have $\hat{T}^{(2)}_x X\cap \hat{T}_y X=\{0\}$,
  then $\textnormal{codim}_{X^*} \sigma_2(X)^*=1$.  In particular $X^*$ is not normal.
\end{prop}
In $\S\ref{homogeneous}$ we apply Proposition \ref{crit} to homogeneous rational varieties. In particular we obtain the following criteria on normality of duals of homogeneous
varieties $G/P$ with $G$ a simple Lie group $G$ and $P$ a parabolic subgroup.
Let $R_+$ denote the set of positive roots (for some choice of the ordering of the roots of the Lie algebra $\gG$) and $w_0$ the involution on the dual of the Cartan subalgebra of $\gG$:
\begin{prop}
 Let $G$ be a simple complex Lie group and $V_\lambda$ an irreducible representation. 
Consider $X=\PP(G.v_\lambda)\subset\PP(V_\lambda)$ the projectivization of 
the highest weight orbit. 
If $X^*$ is normal then either $\sigma_2(X)$ is defective (i.e. not of the expected dimension) or there exists $\alpha,\beta,\gamma \in R_+$ such that 

\[\lambda-w_0(\lambda) =\a+\beta+\g \hspace{1cm}(\diamond)\]
\end{prop}
The table of homogeneous varieties satisfying equation $(\diamond)$ is given and the case of homogeneous varieties with defective secant is detailed. 

In $\S\ref{grass}$ Proposition \ref{crit} and an explicit calculation of the second fundamental form of $\sigma_2(G(3,n))$  allow us to prove Theorem \ref{maintheo}. %The proof is divided in three steps, first the case $k\geq 4$ 
We provide in that section geometric interpretations for the orbits in $\PP(\Lambda^3\CC^8)$ and explicitly describe the bijection between orbits in 
$\PP(\Lambda^3\CC^8)$ and  orbits in the dual space. 
The orbits and their Bruhat order are written down in a graphical way in the appendix of the paper.

In $\S\ref{hyper}$ we show how Proposition \ref{crit} can be applied to Veronese embedding ($v_d(X)$ denotes the $d$-th Veronese re-embedding of  $X$) and Segre products of nondegenerate smooth projective varieties:

\begin{theo}\label{veroseg}
 Let $X\subset \PP^n$ and $Y\subset \PP^m$ two smooth nondegenerate projective varieties. Then 
\begin{enumerate}
 \item For $d\geq 2$, $\sigma_2(v_d(X))^*$ is a codimension one subvariety of $v_d(X)^*$ if and only if $(X,d)\neq (\PP^n,2)$.
\item $\sigma_2(X\times Y)^*$ is a codimension one subvariety of $(X\times Y)^*$ when 
$Y\neq \PP^m$ or $\sigma_2(X)$ is not defective.
\end{enumerate}
\end{theo}

As an example we recover the fact that $(\PP^1\times\PP^1\times \PP^1)^*$ is the only one hyperdeterminant to be normal.

\section{Notations and definitions}\label{def}

\subsection{Second fundamental form} We work throughout with algebraic varieties over the field
$\CC$ of complex numbers. In particular we denote by $V$ a complex
vector space of dimension $N+1$ and $X^n\subset \PP(V)=\PP^{N}$ is a
complex projective nondegenerate variety (i.e. not contained in a hyperplane)
of dimension $n$. Given $x$ a smooth point of $X$, we denote by $T_x
X$ the intrinsic tangent space, $\tilde{T}_x X$ the embedded tangent
space, and
 $\tilde{T}^{(2)}_x X$  the second osculating space of $X$ at $x$.  The notation $\hat{X}$
(resp. $\hat{T}_x X$, ...) means we consider the cone over $X$
(resp. the cone over the embedded tangent space, ...). Define
{\em the conormal space} $N^* _x X:=\hat{x}\otimes (V/\hat{T}_x X)^*$. To
avoid unnecessary complications we ignore twists and write $N^* _x X=(V/\hat{T}_x X)^*$. Let $H$
be a hyperplane tangent to $X$ at $x$, i.e. $\tilde{T}_x X \subset H$, we
denote by $L_H$ the linear form on $V$ defining $H$,  $L_H \in N_x ^* X$ and the restriction $L_{H _{|X,x}}=0$ is a
singular polynomial.  Denote by $S^2 T^*_x X$ the space of quadratic forms on
$T_xX$. The quadratic part of $L_{H _{|X,x}}=0$, denoted
by $Q^H$, allows 
us to define a map, 

\begin{center}
\begin{tabular}{cccc}
${II}_{X,x}:$& $N_x^* X$ & $\onto$ & $S^2 T^*_x X$\\

&$L_H$&$\onto$& $ Q^H$
\end{tabular}
\end{center}

 This map is {\em the
second fundamental form}. Its image is a system of quadrics denoted by
$|{II}_{X,x}|$. We have $|II_{X,x}|\simeq N^* _{2, x}X$ where {\em the
second conormal space} is defined by $N^* _{2, x} X=(\hat{T}^{(2)} _x
X/\hat{T}_x X)^*$  (see \cite{IL}).
By abuse of notation we 
write $H\in
N^* _x X$ instead of $L_H \in N^* _x X$. 

We say $x\in X$ is a {\em
  general point} in the sense of the Zariski topology. The locus of
smooth points of $X$ is denoted by $\text{Sm}(X)$ and the locus of
singular points by $\text{Sing}(X)$.

\subsection{Auxiliary varieties} The s-secant variety of a projective variety $X\subset \PP^{m}$ is the variety $\sigma_s(X)$ defined to be 
the Zariski closure of the union of the linear span of $s$-tuples points of $X$
\[\sigma_s(X)=\overline{\bigcup_{x_1,\dots,x_s\in X} \PP^{s-1}_{x_1,\dots,x_s}}\]
where $\PP^{s-1}_{x_1,\dots,x_s}$ is a projective space of dimension $s-1$ passing through $x_1,\dots, x_s$.
The dimension of $\sigma_s(X)$ is often calculated from the famous Terracini's Lemma \cite{Z2},
\begin{theo}{\bf [Terracini's Lemma]}: Let $x_1,\dots, x_s$ be a general collection of points of $X$ and let $z$
 be a general point in $\PP^{s-1} _{x_1,\dots,x_s}$. Then the 
tangent space to $\sigma_s(X)$ at $z$ is given by \[\tilde{T}_z \sigma_s(X)=<\tilde{T}_{x_1}X,\dots,\tilde{T}_{x_s} X>\]
where $<\tilde{T}_{x_1}X,\dots,\tilde{T}_{x_s} X>$ denotes the projective span.
\end{theo}
\begin{rem}
 It is clear from Terraccini's Lemma that given  $H$ a smooth point of $\sigma_s(X)^*$, i.e. $H$ is a <<general>> hyperplane tangent 
to $\sigma_s(X)$, then $H$ is tangent to $X$ at $s$ points. In other words Terracini's Lemma implies $\sigma_s(X)^*\subset X^*$.
\end{rem}

Let $X$ and $Y$ be two projective varieties and let $\PP^1 _{xy}$ denote the projective line containing $x\in X$ and $y\in Y$. 
The join of $X$ and $Y$ is  the Zariski closure of the lines joining $X$ and $Y$:
\[J(X,Y)=\overline{\bigcup_{x\in X, y\in Y, x\neq y} \PP^{1}_{xy}}\]

In particular $J(X,X)=\sigma_2(X)$. 

Assume $Y\subset X$ and let $T^{\star}_{X,Y,y_0}$ denote the union of $\PP^1 _*$'s where $\PP^1 _*$ is 
the limit of $\PP^1_{xy}$ with $x\in X$, $y\in Y$ and $x,y\to y_0\in Y$. 
The union of the $T^{\star}_{X,Y,y_0}$ is called the variety of  relative tangent stars of $X$ with respect to $Y$ (see \cite{Z2}):
\[T(Y,X)=\bigcup_{y \in Y}T^{\star}_{X,Y,y} \]
If $Y=X$, the variety $T(X,X)$, also denoted by $\tau(X)$, is the usual tangential variety.

\setcounter{prop}{0}
\section{Dimension of $\sigma_2(X)^*$}~\label{comp}

In this section we give a sufficient condition for $\sigma_2(X)^*$ to be of maximal dimension in $X^*$.

\begin{lemm}~\label{sec}
Let $X\subset \PP(V)$ be a smooth projective
variety of dimension $n$. Suppose for a general pair of point $(x,y)\in X\times X$ there exists $H$ with the
following properties,
\[\left\{\begin{array}{l}
\tilde{T}_x X\subset H, \tilde{T}_y X \subset H\\
          \textnormal{rank}(II_{X,x}(H))=n \text{ and } \textnormal{rank}(II_{X,y}(H))=n
         \end{array}\right.\]

Then $\textnormal{codim}_{X^*}\sigma_2(X)^*=1$.
\end{lemm}

\proof 
The Katz dimension formula \cite{Katz} gives the dimension of the dual of a projective variety 
from the rank of a generic quadric in the image of the second 
fundamental form. More precisely let
$Z\subset \PP(V)$ be a projective variety of codimension $a$ and $z\in
Z$ is a general point. Assume $r$ is the rank
of a generic quadric in $|{II}_{Z,z}|$, then
$\text{dim}(Z^*)=r+a-1$. In particular when generic quadrics in $|II_{Z,z}|$ are not of maximal rank, the dimension of $Z^*$ is less than expected, i.e. $Z^*$ is not 
a hypersurface and the tangent hyperplanes are tangent to $Z$ along a (at least) one dimensional subset of $Z$.

Let us consider
$Z=\sigma_2(X)$. By the Terracini lemma, if $(x,y)$ is a general pair of points
of $X\times X$, then the tangent space $\tilde{T}_z\sigma_2(X)$ is
constant along the line $\PP^1_{xy}$ (by abuse of notation we will write $z=x+y$ the point $z$ on $\PP_{xy}^1$). Therefore a quadric in
$|{II}_{\sigma_2(X),z}|$ has always a degenerate direction and its rank is
bounded by $2n$. Suppose there exists $H$
with the hypothesis of the lemma. Then we claim that
$\text{rank}({II}_{\sigma_2(X),z}(H))=2n$. If not there exists a curve
$z(t) \not\subset\PP^1_{xy}$ such that $H$ is tangent to $\sigma_2(X)$ along
$z(t)$. But $z(t)=x(t)+y(t)$ and we can suppose $x(t) \not \equiv x$. Then the hyperplane $H$ is tangent to $X$ along $x(t)$, but it contradicts
the assumption $
\textnormal{rank}({II}_{X,x}(H))=n$. $\Box$\\

\begin{rem}
The hypothesis  on the
rank of $II_{X,x}(H)$ implies that $X^*$ is a hypersurface.
\end{rem}

We now state our criteria to have $\text{codim}_{X^*} \sigma_2(X)^*=1$:

\begin{prop}~\label{crit}
Let $X\subset \PP(V)$ be a smooth projective variety of dimension $n$.
Assume $X^*$ is a hypersurface. Suppose for a general pair of point $(x,y)\in X\times
X$ we have $\hat{T}^{(2)}_x X\cap \hat{T}_y X=\{0\}$,
  then $\textnormal{codim}_{X^*} \sigma_2(X)^*=1$.  In particular $X^*$ is not normal.
\end{prop}

\begin{rem}
In  Proposition \ref{crit}, a consequence of the hypothesis
$\hat{T}^{(2)}_x X\cap \hat{T}_y X=\{0\}$ is that $\sigma_2(X)$
is of maximal dimension (we also say nondefective). %Terraccini
\end{rem}

\proof 
Let $z\in \PP^1 _{xy}\subset \sigma_2(X)$ be a general point of the 2-secant variety. Let
us consider the maps:

$$r: N^* _x X \twoheadrightarrow N^* _{2,x} X=(\hat{T}_x ^{2} X/\hat{T}_x X)^*$$

$$i:N_z ^* \sigma_2(X)= (V/<\hat{T}_{x} X+\hat{T}_{y} X>)^* \hookrightarrow  N_x^* X.$$

The assumption $\hat{T}^{(2)} _{x} X\cap \hat{T}_{y} X
=\{0\}$ says that for any $L\in N^*_{2,x} X$ one can find a hyperplane
$H\in N^*
_x X$ such that its restriction is $L=r(H)$, and
$\hat{T}_{y}X\subset H$ (i.e. $H$ is obtained, by the map $i$ from a
hyperplane of $N_z ^* \sigma_2(X)$, we write $H\in N_z ^* \sigma_2(X)$). The dual variety $X^*$ is a hypersurface by
hypothesis, thus we can choose $L$
such that $II_{X,x}(H)=Q^L$ is of full rank. One obtains a hyperplane $H \in
N^*_z \sigma_2(X)$ such that the quadric $Q^H=Q^{L}\in |II_{X,x}|$
is of rank $n$. The same construction works for $y$ and one obtains a
hyperplane $H'\in N^* _z\sigma_2(X)$ such that $II_{X,y}(H')=Q^{L'}$ is
of rank $n$. We now consider the line $\PP^1
_{HH'}$. That line can be seen either as a line in the
projectivized conormal space
of $X$ at $x$  ($\PP^1 _{HH'} \subset
\PP(N^*_{x }X)$) or in the projectivized conormal space of $X$ at $y$ ($\PP^1 _{HH'} \subset \PP(N^*_{y} X)$). In each case
 there is only a finite number  of points on $\PP^{1} _{HH'}$ such that the
corresponding quadrics at $x$ and $y$ are not of full rank (the
quadrics corresponding to $H$ and
$H'$ being of full rank, the line can not be contained in the
subvariety of degenerate quadrics). In other words there exists $H''\in
\PP^1_{HH'}$ (there exists an infinity of such) such that $II_{X,x}(H'')$ and $II_{X,y}(H'')$ are of rank
$n$. The lemma \ref{sec} implies $\sigma_2(X)^*$ has codimension $1$ in $X^*$. $\Box$

\section{Application: A criteria of Normality for $(G/P)^*$}\label{homogeneous}
We apply Proposition \ref{crit} when $X=G/P$ is a rational homogeneous variety. One obtains a general criteria which is a necessary condition for $(G/P)^*$ to be normal
and we list the homogeneous varieties which satisfy the criteria.
Let $G$ be a complex simple Lie group and $P$ a parabolic subgroup. % corresponding to a subset $S$ of nodes of the  Dynkin diagram. 
The homogeneous space $G/P$ has 
a homogeneous embedding in an irreducible representation $V_\lambda$ of $G$ (representation of highest weight $\lambda$) where $\lambda=\sum_{i}a_i\omega_i$
 with $\omega_i$ the $i$-th fundamental weight and $a_i\in \NN$. The embedding $X=G/P\subset \PP(V_\lambda)$ is the projectivization of the highest weight orbit,
 i.e. $X=\PP(G.v_\lambda)\subset\PP(V_\lambda)$. For $\lambda=\omega_i$ we denote by $P_i$ the corresponding parabolic subgroup.% and $G/P_i$ the corresponding homogeneous variety.

\begin{ex}
 \begin{enumerate}
  \item Let $G=SL_n$, which acts on $V=\CC^n$, and $1\leq a_1<a_2<\dots<a_p\leq n-1$. The Lie group $G$ acts on  $W=\Lambda^{a_1} V\otimes \Lambda^{a_2}V\otimes \dots \otimes\Lambda^{a_p}V$, and the 
highest weight orbit in $\PP(W)$ is the flag variety $\mathbb{F}_{a_1,\dots,a_p}(V)$, i.e. the variety of (partial) flag $0\subset E_1\subset\dots\subset E_p\subset V$ with $E_i$ linear space 
of $V$ such that $\text{dim}(E_i)=i$. In particular the variety $\mathbb{F}_k(V)$ is the Grassmannian of $k$-planes in $V$. The variety of complete flag $\mathbb{F}_{1,\dots,n-1}(V)$ is obtained with 
$\lambda=\omega_1+\dots+\omega_{n-1}$ and $P=B$ is the Borel subgroup of $G$ see \cite{F-H}.
\item Let $G=SO_{n}$ acting on $V=\CC^n$, $V$ equipped with a nondegenerate quadratic form $Q$, $W=\Lambda^k V$ and $\lambda=\omega_k$. 
The corresponding highest weight orbit $G/P_{k}\subset \PP(W)$ is the variety $X=G_Q(k,n)$ 
the Grassmannian of isotropic $k$-planes, i.e. $X=G_Q(k,n):=\{E\in G(k,n), Q(v,w)=0\ \forall\ v,w\in E\}$. 
For $k=m-1$ and $n=2m-1$ (resp. $k=m$ and $n=2m$) the variety $G_Q(m-1,2m-1)$ (resp. $G_Q(m,2m)$) has two isomorphic components. 
The components are called Spinor varieties $\SS_m$ and can be obtained as the highest weight orbit of the spinor representation of the group $Spin_{2m-1}$ of type $B_{m-1}$ with highest weight $\omega_{m-1}$
 (resp. $Spin_{2m}$ of type $D_m$ with highest weight $\omega_m$).

\item Let $G=Sp_{2n}$ acting on $V=\CC^{2n}$, $V$ equipped with a nondegenerate symplectic form $\omega$, and $W=\Lambda^k W$. 
The variety $X=G/P_k\subset \PP(W)$ is the Grassmannian of isotropic $k$-planes for $\omega$, $X=G_\omega(k,n):=\{E\in G(k,n),\omega(v,w)=0\ \forall\ v,w\in E\}$.
%\item The Spinor variety  
\item Let $G=E_6$ and $\lambda=\omega_1$. The homogeneous variety $\mathcal{E}_6=E_6/P_1\subset\PP(V_{\omega_1})$, called the Severi variety of type $E_6$, 
can be identified with the Cayley projective plane $\mathbb{O}\PP^2$ embedded in the Jordan algebra of
the $3\times 3$ $\mathbb{O}$-Hermitan symmetric matrices, see \cite{LM2,Z2}.
 \end{enumerate}

\end{ex}

For $X=\PP(G.v_\lambda)\subset \PP(V_\lambda)$ a general pair of points can be chosen to be $(v_\lambda,v_{\mu})$ with $\mu$ the lowest weight of 
the representation.
The representation-theoretic interpretation of the osculating spaces of $X$ are given in \cite{LM2}: consider $\gG$ the simple Lie algebra of $G$ and 
$\gG^{(2)}$, the second term in the natural filtration of the universal Lie algebra, i.e. 
$\gG^{(2)}=\gG\otimes\gG/\{x\otimes y-y\otimes x-[x,y]| x,y\in \gG\}$, then the tangent and second osculating spaces at $v_\lambda$ and $v_{\mu}$ are given by
\[\hat{T}_{v_\lambda}=\gG. v_\lambda \text{ and }\hat{T}^{(2)}_{v_{\mu}}X=\gG^{(2)}.v_\mu\]

Moreover if we denote by $R_+$ the positive roots of $\gG$ and by $V_\rho$ the eigenspace corresponding to the weight $\rho$ we have 
(the root spaces $\gG_\alpha$ of $\gG$ act on the eigenspaces $V_\rho$ by <<translation>> see \cite{F-H}):

\[\gG. v_{\lambda}\subset v_{\lambda}\oplus(\oplus_{\gamma\in R_+} V_{\lambda-\gamma}) \text{ and } \gG ^{(2)}. v_{\mu}\subset v_{\mu}\oplus(\oplus_{\alpha\in R_+}V_{\mu+\alpha})\oplus(\oplus_{\alpha,\beta\in R_+} V_{\mu+\alpha+\beta}) \] 

The condition $\hat{T}_{v_{\lambda}}X\cap \hat{T}^{(2)}_{v_\mu}X=\{0\}$ is satisfied
 when $\lambda-\gamma\neq \mu+\alpha+\beta$ for all $\alpha,\beta,\gamma \in R_+$ and $\lambda-\gamma\neq \mu+\alpha$ for all $\alpha,\gamma \in R_+$ 
(this corresponds to $\hat{T}_{v_{\lambda}}X\cap \hat{T}_{v_{\mu}}X=\{0\}$, i.e. $\sigma_2(X)$ is nondefective). 
Denote by $w_0$ the involution on the dual of the Cartan subalgebra of
$\gG$, which transforms $R_+$ into $R_-$, i.e such that $w_0(\lambda)=\mu$ then a consequence of Proposition \ref{crit} is the following general statement:

\begin{prop}
 Let $G$ be a simple complex Lie group and $V_\lambda$ an irreducible representation. 
Consider $X=\PP(G.v_\lambda)\subset\PP(V_\lambda)$ the projectivization of 
the highest weight orbit. 
If $X^*$ is normal then either $\sigma_2(X)$ is defective or there exists $\alpha,\beta,\gamma \in R_+$ such that 

\[\lambda-w_0(\lambda) =\a+\beta+\g \hspace{1cm}(\diamond)\]
\end{prop}

In Table~\ref{table} we list the homogenous varieties $G/P$ which satisfy $(\diamond)$, i.e. such that their duals $(G/P)^*$ are potentially normal.
\begin{ex}
 As an example we solve equation $(\diamond)$ when $G=F_4$ (the proof follows the same steps for the other types). We use the notation of \cite{Bou} and
 denote by $W=\RR^4$ the real vector space spanned by the root lattice of the Lie algebra of type $F_4$ with orthogonal basis  $(\epsilon_i)_{1\leq i\leq 4}$. The positive
roots are $\epsilon_i,\epsilon_i\pm\epsilon_j$ ($i<j$) and $\frac{1}{2}(\epsilon_1\pm\epsilon_2\pm\epsilon_3\pm\epsilon_4)$. The fundamental weights are $\omega_1=\epsilon_1+\epsilon_2$, 
$\omega_2=2\epsilon_1+\epsilon_2$, $\omega_3=\frac{1}{2}(3\epsilon_1+\epsilon_2+\epsilon_3+\epsilon_4)$ and $\omega_4=\epsilon_1$.
 We consider the following norm on $W$:
for all $v\in W$ such that $v=\sum_{i=1}^4 a_i\epsilon_i$ we define $||v||=\sum_{i=1}^4 |a_i|$. In particular if $\alpha$ is a positive root 
we have $||\alpha||=1 \text{ or }2$ and the norm of the sum of three positive roots $||\alpha+\beta+\gamma||$
is equal to $3, 4,5$ or $6$.
On the other hand the involution $w_0$ on $W$ gives $w_0(\epsilon_i)=-\epsilon_i$. Thus for any fundamental weight 
we have $||\omega_i-w_0(\omega_i)||=||2\omega_i||$ which can be equal to $2$ ($i=4$), $4$ ($i=1$), or $6$ ($i=2,3$). 
Let $\lambda$ be a highest weight  for the Lie algebra $F_4$, then $\lambda=\sum_{i=1}^4 a_i\omega_i$ with $a_i\geq 0$ and
 $||\lambda-w_0(\lambda)||=\sum_{i=1}^4 |a_i|||\omega_i-w_0(\omega_i)||$. The restrictions on the possible values of $||\alpha+\beta+\gamma||$ 
and $||\omega_i-w_0(\omega_i)||$ 
lead to the following list of weights $\lambda$ which can be solution of $(\diamond)$: $\lambda=3\omega_4,2\omega_4,\omega_3,\omega_2,\omega_1,\omega_1+\omega_4$.
Among those candidates one checks that only $\omega_3-w_0(\omega_3)$ is the sum of three positive roots. 
Thus $\lambda=\omega_3$ is the only solution for the Lie group of type $F_4$.$\Box$
\end{ex}

\begin{table}[h]
\caption{Homogeneous varieties which satisfy $(\diamond)$}\label{table}
\begin{tabular}{cccc}
   \hline
Type & highest weight & $X$ \\
\hline
$A_n$ & $3\omega_1, 3\omega_n$ & $v_3(\PP^n)$\\
        & $\omega_1+\omega_2, \omega_{n-1}+\omega_n$ & $\mathbb{F}_{1,2}(\CC^n)$\\
       & $\omega_1+\omega_{n-1},\omega_2+\omega_n$ & $\mathbb{F}_{1,n-1}(\CC^n)$\\
       & $\omega_3$ &   $G(3,n)$\\
\hline
$B_5$&$\omega_5$ &  $\SS_6$\\
 $B_6$     &$\omega_6$ &  $\SS_7$\\
 $B_n$     &$\omega_3$   & $G_Q(3,2n+1)$\\
\hline
$C_n$ & $3\omega_1$ &  $v_3(\PP^{2n-1})$\\
      & $\omega_1+\omega_2$              &$\mathbb{F}_{1,2,\omega}(\CC^{2n})\text{ (isotropic flag variety)}$\\
        &$\omega_3$ & $G_\omega (3,2n)$ \\
\hline
$D_6$ & $\omega_6$ &  $\SS_6$\\
$D_7$ & $\omega_7$ &  $\SS_7$\\
$D_n$ & $\omega_3$ & $G_Q(3,2n)$\\
\hline
$E_6$ &$\omega_3,\omega_5$ &  $E_6/P_3$\\
\hline
$E_7$ &$\omega_2$ & $E_7/P_2$\\
      &$\omega_7$ &  $E_7/P_7$\\
\hline
$F_4$ & $\omega_3$ & $F_4/P_3$\\
\hline
  \end{tabular}
\end{table}

We now discuss the normality of the dual of $(G/P)$ when $(G/P)$ has a defective secant variety. 
The homogeneous varieties with defective secant are (see \cite{Ka}): the smooth quadric hypersurface  $\QQ^n$, the $5$-Spinor variety $\SS_5$, 
the Scorza varieties ($v_2(\PP^n)$, $\PP^n\times\PP^n$, $G(2,n)$, $\mathcal{E}_6$), 
the general hyperplane section of $\mathcal{E}_6$  and the adjoint varieties (highest weight orbit for the action of $G$ on $\PP(\gG)$) . 

\begin{enumerate} 
 \item The varieties $\QQ^n, \SS_5$ are smooth self-dual varieties and therefore normal.
\item The varieties $\PP^n\times\PP^n$, $G(2,n)$ have normal duals as it has been recalled in the introduction. The same is true for $v_2(\PP^n)$ as it will be stated in $\S\ref{hyper}$.
\item A description of the varieties  $\mathcal{E}_6=E_6/P_1$ and $F_4/P_4=\mathcal{E}_6\cap H$ (general hyperplane section of $\mathcal{E}_6$) and their duals can be found in \cite{Z2} Chapter III. The dual of the Severi variety $\mathcal{E}_6$ is normal
and the dual of $\mathcal{E}_6\cap H$ is not.
 %\item The variety $F_4/P_1=\mathcal{E}_6\cap H$ is not normal.
\item The normality of the duals of the adjoint varieties  can be solved using results of \cite{Kn}. In that paper
F. Knop studied the hyperplane sections of the adjoint varieties. Because a normal variety is regular in codimension $1$, 
a  dual hypersurface $X^*$ is normal if and only if it parametrizes 
hyperplane sections of $X$ with either a unique quadratic singularity (hyperplanes which correspond to smooth points of $X^*$) or with  nonisolated singularities 
(hyperplanes in $\text{Sing}(X^*)$). Following  Knop's Theorem only the adjoint varieties for the Lie
groups $Sp_n$ (i.e. $X=v_2(\PP^n)$ see $\S 6$) and $G_2$ (denoted by $X=G_2/P_2$) have singular hyperplane sections with either a unique quadratic singularity or nonisolated singularities. Thus the only adjoint varieties with normal duals are $v_2(\PP^n)$ and $G_2/P_2$.
\end{enumerate}

\begin{rem}
 Some varieties of Table \ref{table} have been studied in details and we know according to \cite{LM1} that the varieties  $G_\omega(3,6)$, $G(3,6)$, $\SS_6$, $E_7/P_7$
have normal duals.
\end{rem}

\section{Normality of the duals of Grassmannians}\label{grass}
In this section we prove Theorem \ref{maintheo} in three steps. First we apply Proposition \ref{crit} 
(we recover without reference to roots and weights the result of $\S\ref{homogeneous}$ in the case of the 
Grassmannians). Then we study in details the case of $G(3,n)$ with $n\geq 9$ where Proposition \ref{crit} does not apply directly. 
The remaing cases correspond to the action of $SL_n$ on $\Lambda^3\CC^n$ with finitely many orbits ($n=6,7,8$). 
In $\S\ref{hyper}$ we will recover the result of \cite{WZ} on normality of hyperdeterminants following the same three steps.
\subsection{The varieties $G(k,n)$ with $k\geq 4$}
Proposition \ref{crit} implies $G(k,n)^*$ is not normal for $k\geq 4$:
\begin{prop}
 Let $X=G(k,n)\subset \PP(V)$, with $k\geq 3$. Given a general pair of points $(x,y)\in G(k,n)\times G(k,n)$ we have
\[\hat{T}_x^{(2)}G(k,n)\cap \hat{T}_y G(k,n)\neq \{0\}\Leftrightarrow k=3\] 
\end{prop}

\proof Consider $E$ and $E'$ two transverse $k$-planes in $V=\CC^n$, i.e. $([E],[E'])$ is a general pair of point in $G(k,n)\times G(k,n)$. The tangent and second osculating spaces
at $E$ and $E'$ are:
\[\hat{T}_E G(k,n)=\Lambda^{k-1}E\Lambda V \text{ and } \hat{T}^{(2)}_{E'}G(k,n)=\Lambda^{k-2} E'\Lambda(\Lambda^2 V)\]

It is clear that $\hat{T}_E G(k,n)\cap \hat{T}^{(2)}_{E'}G(k,n)=\{0\}$ for all $k$ such that $k-2\geq 2$.$\Box$

\begin{coro}\label{k4}
 If $k\geq 4$ then $G(k,n)^*$ is singular in codimension one.
\end{coro}

\subsection{The varieties $G(3,n)$ with $n\geq 9$}

In the case of Grassmannians of $3$-planes Proposition \ref{crit} does not allow us to conclude. However the proof of Proposition \ref{crit} is based on the existence of
a hyperplane $H\in N_{x}^*X$ such that $r(H)$ is of maximal rank and  $H\supset T_y X$. We now prove the existence of such a hyperplane for $k=3$ and $n\geq 9$.
Let $e_1,\dots,e_n$ a basis of $V=\CC^n$. Using Pl\"ucker embedding we denote a general pair of points $(x,y)$  by
 $x=[e_1\wedge e_2\wedge e_3]$ and $y=[e_4\wedge e_5\wedge e_6]$  ($U=<e_1, e_2,e_3>$ and $U'=<e_4, e_5, e_6>$ be the corresponding $3$-planes in $V=\CC^n$). A direct calculation shows that 
\[\hat{T}_x ^{(2)} G(3,n)\cap \hat{T}_y G(3,n)=U\Lambda(\Lambda^2U')\]
Thus $H\in N_x ^* G(3,n)$ is tangent to $G(3,n)$ at $y$ if and only if 
\[r(H) \in (\hat{T}_x ^{(2)}G(3,n)/(\hat{T}_x G(3,n)+\hat{T}_x^{(2)} G(3,n)\cap \hat{T}_y G(3,n)))^*\]
i.e. \[r(H)\in (U\Lambda(\Lambda^2V)/(\Lambda^2 U\Lambda V+U\Lambda(\Lambda^2U')))^* \]
Given such $H$ can we have $II_{G(3,n),e_1\wedge e_2\wedge e_3}(H)$ is a quadric of full rank?
To answer that question one needs to compute $II_{G(3,n),e_1\wedge e_2\wedge e_3}$. 
This can be done using moving frames techniques (see \cite{IL} page 100):
\[II_{G(3,n),e_1\wedge e_2\wedge e_3}= \sum_{s<t}(\omega_2 ^s \omega_3 ^t-\omega_2 ^t \omega_3 ^s)e_1\wedge e_s\wedge e_t+(\omega_1 ^s \omega_3 ^t-\omega_1 ^t \omega_3 ^s)e_s\wedge e_2\wedge e_t+
(\omega_2 ^s \omega_1 ^t-\omega_2 ^t \omega_1 ^s)e_1\wedge e_2\wedge e_t\]
where $\omega=(\omega_i ^s)$ is the Maurer-Cartan form for the $GL(\CC^n)$-frame bundle with indices $1\leq i \leq 3$ and $4\leq s\leq n$ .
Considering $\{\omega_i^s\}$ as a basis of $T_{e_1\wedge e_2\wedge e_3}^* G(3,n)$, then for any $H\in N_x^* G(3,n)$ the quadric 
$II_{G(3,n),e_1\wedge e_2\wedge e_3}(H)\in S^2 T^*_x G(3,n)$ is of type \[Q^H=\begin{pmatrix}
                                                                                                                               0 & A & B \\
^t A & 0 & C\\
^tB & ^t C & 0
                                                                                                                              \end{pmatrix} \hspace{1cm}(\star)\]
with $A,B,C$ being skew symmetric matrices of size $(n-3)\times (n-3)$.

The condition $r(H)\in (U\Lambda(\Lambda^2V)/(\Lambda^2 U\Lambda V+U\Lambda(\Lambda^2U')))^* $  
 is equivalent to the fact that there are no terms of type
$\omega_1^s\omega_2^t-\omega_2^s\omega_1^t$, $\omega_1^s\omega_3^t-\omega_3^s\omega_1^t$ and $\omega_2^s\omega_3^t-\omega_3^s\omega_2^t$ with $4\leq s<t\leq 6$ in $II(H)$. 
In other words the matrices $A,B,C$ are skew-symmetric and of type:
\[\left(\begin{array}{ccc|ccccc}
           0 & 0& 0& & & & & \\
           0 &0 & 0& & &*& &\\
           0 & 0 &0 & & &&&\\
\hline
             & & &  & & & & \\
    & & &  & & & & \\
&* & &  & & *& & \\
& & &  & & & & \\
          \end{array}\right)\hspace{1cm} (\star\star)\]
One needs to determine if one can build a symmetric matrix $Q^H$ of maximal rank of type $(\star)$ with the condition $(\star\star)$ on the blocks. If such a quadric exists then there exists $H$ such that $H$ is tangent to $X$ at $x$ and $y$ and $Q^H$ is of maximal rank. 
\begin{lemm}\label{calc}
 For $n\geq 9$ such a quadric exists.
\end{lemm}

\proof By induction
\begin{itemize}

\item From $n$ to $n+3$: suppose $Q$ is a quadric satisfying $(\star)$ and $(\star\star)$. We consider the basis $\{\omega_1 ^s,\omega_2 ^s,\omega_3^s,\omega_1 ^{n+i},\omega_2 ^{n+i},\omega_3^{n+i}\}$ 
with $4\leq s\leq n$ and $1\leq i\leq 3$. Then the following quadric of size $(3(n+3))^2$ satisfies $(\star)$ and $(\star\star)$ : 
$\left(\begin{array}{ccc|c}
 & & &0 \\
 &Q& &\vdots \\
 &   & &0 \\
\hline
 0 &\dots &0 &q \\

                                                                                                       \end{array}\right)$
 where $q$ is a symmetric matrix of size $9\times 9$ and rank $9$ satisfying $(\star)$.
\item For $n=9,10,11$ we give explicit examples of symmetric matrices satisfying $(\star)$, $(\star\star)$:
\begin{enumerate}
 \item $n=9$ we consider the blocks:

\[A=\begin{pmatrix}
0 & 0 & 0 & -1 & 0 & 0  \\
 0 & 0 & 0 & 0 & -1 & 0 \\ 
 0 & 0 & 0 & 0 & 0 & -1 \\
 1 & 0 & 0 & 0 & 0 & 0 \\
 0 & 1 & 0 & 0 & 0 & 0 \\ 
 0 & 0 & 1 & 0 & 0 & 0\end{pmatrix}, B=\begin{pmatrix}
0 & 0 & 0 & -1 & 0 & 0\\
0 & 0 & 0 & -1 & -1 & 0\\
 0 & 0 & 0 & 0 & -1 & -1\\
1 & 1 & 0 & 0 & 0 & 0\\
0 & 1 & 1 & 0 & 0 & 0\\
0 & 0 & 1 & 0 & 0 & 0
\end{pmatrix}, C=\left(\begin{array}{cccccc}
 0 & 0 & 0 & 0 & -1 & 0\\
 0 & 0 & 0 & 0 & 0 & 1\\
 0 & 0 & 0 & -1 & 0 & 0\\
 0 & 0 & 1 & 0 & 0 & 0\\
 1 & 0 & 0 & 0 & 0 & 0\\
 0 & -1 & 0 & 0 & 0 & 0
  \end{array}\right)\]

\item For $n=10$,
\[A=\begin{pmatrix}
    0 & 0 & 0 & 1 & 0 & 0 & 1 \\
0 & 0 & 0 & 0 & 1 & 0 & 1 \\
 0 & 0 & 0 & 0 & 0 & 1 & 1 \\
  -1 & 0 & 0 & 0 & 0 & 0 & 1 \\
 0 & -1 & 0 & 0 & 0 & 0 & 1 \\
 0 & 0 & -1 & 0 & 0 & 0 & 1 \\
 -1 & -1 & -1 & -1 & -1 & -1 & 0  
    \end{pmatrix}, B=\begin{pmatrix}
0 & 0 & 0 & 1 & 1 & 0 & 1\\
 0 & 0 & 0 & 0 & 1 & 1 & 1\\   
 0 & 0 & 0 & 0 & 0 & 1 & 1\\
-1 & 0 & 0 & 0 & 0 & 0 & 1\\ 
-1 & -1 & 0 & 0 & 0 & 0 & 1\\ 
0 & -1 & -1 & 0 & 0 & 0 & 1\\
-1 & -1 & -1 & -1 & -1 & -1 & 0
    \end{pmatrix}\]

\[C=\begin{pmatrix}
 0 & 0 & 0 & 0 & 1 & 0 & 1\\ 
 0 & 0 & 0 & 0 & 0 & -1 & 2\\ 
0 & 0 & 0 & 1 & 0 & 0 & 3\\
 0 & 0 & -1 & 0 & 0 & 0 & 4\\
-1 & 0 & 0 & 0 & 0 & 0 & 5\\
  0 & 1 & 0 & 0 & 0 & 0 & 6\\ 
-1 & -2 & -3 & -4 & -5 & -6 & 0

    \end{pmatrix}\]

\item For $n=11$

\[A=\begin{pmatrix}
 0 & 0 & 0 & 0 & -1 & 0 & 0 & 0 \\ 
 0 & 0 & 0 & 0 & 0 & -1 & 0 & 0 \\
 0 & 0 & 0 & 0 & 0 & 0 & -1 & 0 \\
 0 & 0 & 0 & 0 & 0 & 0 & 0 & -1 \\ 
 1 & 0 & 0 & 0 & 0 & 0 & 0 & 0 \\
 0 & 1 & 0 & 0 & 0 & 0 & 0 & 0 \\
 0 & 0 & 1 & 0 & 0 & 0 & 0 & 0 \\
  0 & 0 & 0 & 1 & 0 & 0 & 0 & 0 

  \end{pmatrix}, B=\begin{pmatrix}
     0 & 0 & 0 & 0 & -1 & 0 & 0 & 0\\
0 & 0 & 0 & 0 & -1 & -1 & 0 & 0\\
0 & 0 & 0 & 0 & 0 & -1 & -1 & 0\\
0 & 0 & 0 & 0 & 0 & 0 & -1 & -1\\
1 & 1 & 0 & 0 & 0 & 0 & 0 & 0\\
0 & 1 & 1 & 0 & 0 & 0 & 0 & 0\\
0 & 0 & 1 & 1 & 0 & 0 & 0 & 0\\
0 & 0 & 0 & 1 & 0 & 0 & 0 & 0
    \end{pmatrix}\]

\[C=\begin{pmatrix}
0 & 0 & 0 & 0 & 0 & -1 & 0 & 0\\
 0 & 0 & 0 & 0 & 0 & 0 & -2 & 0\\
 0 & 0 & 0 & 0 & 0 & 0 & 0 & -3\\
 0 & 0 & 0 & 0 & 0 & 0 & 0 & 0\\
 0 & 0 & 0 & 0 & 0 & 0 & 0 & 0\\
 1 & 0 & 0 & 0 & 0 & 0 & 0 & 0\\
 0 & 2 & 0 & 0 & 0 & 0 & 0 & 0\\
 0 & 0 & 3 & 0 & 0 & 0 & 0 & 0\\
    \end{pmatrix}\]

\end{enumerate}

\end{itemize}
\begin{coro}\label{k9}
 If $n\geq 9$ then $G(3,n)^*$ is singular in codimension $1$.
\end{coro}

\subsection{The varieties $G(3,6)$, $G(3,7)$ and $G(3,8)$}\label{orbits}
To complete the proof of Theorem \ref{maintheo} we prove that $G(3,6)^*$, $G(3,7)^*$ and $G(3,8)^*$ are regular in codimension one (the cases $G(3,5)$ and $G(3,4)$ follow from $G(3,5)\simeq G(2,5)$ and $G(3,4)\simeq \PP^3$). 
The classification of orbits for the action of $SL_{n}$ on $\Lambda^3 \CC^n$ is known for $n\leq 8$ (there is a classification for $n=9$ but the number of orbits is not finite) see \cite{Djo,Gu}.
For $n=6,7$ the geometric nature of the orbits has been investigated in various papers (\cite{Do,LM1} for $n=6$ and \cite{AOP} for $n=7$). 
However for $n=8$ we did not find in the literature a geometric approach for the orbits decomposition. We take advantage of the 
present paper to put together what is known for $n=6,7$ and provide geometric descriptions for $n=8$. In particular we describe 
the duality between the orbits and answer in this particular case a question of E. A. Tevelev on group actions with finitely many orbits (see $\S2.2$ of \cite{Te} on the Pyasetskii Pairing). To simplify the notation we will write $e_{ijk}$ for $e_i\wedge e_j\wedge e_k$.
The following tables give for each orbit a representative, the dimension and the variety corresponding to the closure of the orbit.
A direct consequence of the tables  is:
\begin{coro}\label{k8}
 The varieties $G(3,6)^*,G(3,7)^*$ and $G(3,8)^*$ are regular in codimension one.
\end{coro}
\begin{rem}
 The corollaries \ref{k4}, \ref{k9} and \ref{k8} prove Theorem \ref{maintheo}.
\end{rem}

\begin{table}[h]
\caption{$SL_6$-orbits in $\PP(\Lambda^3\CC^6)$}\label{t6}
\begin{tabular}{|c|c|c|}
  \hline
\text{Representative} & \text{Dimension} & \text{Geometric interpretation} \\
\hline
$e_{123}$ & $9$& $G(3,6)$\\
\hline
$e_{123}+e_{345}$ & $14$ & $\text{Sing}(G(3,6)^*)\simeq\sigma_{2,+}(G(3,6))$\\
\hline
$e_{123}+e_{345}+e_{156}$& $18$ & $\tau(G(3,6))\simeq G(3,6)^*$\\
\hline
$e_{123}+e_{456}$ & $19$ & $\PP^{19}$\\
\hline
\end{tabular}
\end{table}

\begin{table}[h]
\caption{$SL_7$-orbits in $\PP(\Lambda^3\CC^7)$}\label{t7} 
\begin{tabular}{|c|c|c|}
  \hline
\text{Representative} & \text{Dimension} & \text{Geometric interpretation} \\
\hline
$e_{123}$ & $12$ & $G(3,7)$\\
\hline
$e_{123}+e_{147}$ & $19$ & $\sigma_{2,+}(G(3,7))$\\
\hline
$e_{456}+e_{147}+e_{257}$ & $24$ & $\tau(G(3,7))$\\
\hline
$e_{123}+e_{456}$ & $25$ & $\sigma_2(G(3,7))$\\
\hline
$e_{147}+e_{257}+e_{367}$ & $20$ & $\sigma_2(G(3,7))^*$\\
\hline
$e_{456}+e_{147}+e_{257}+e_{367}$ & $27$ & $\tau(G(3,7))^*$\\
\hline
$e_{123}+e_{456}+e_{147}$ & $30$ & $\text{Sing}(\sigma_3(G(3,7))\simeq\sigma_{2,+}(G(3,7))^*$\\
\hline
$e_{123}+e_{456}+e_{147}+e_{257}+e_{367}$ & $33$ & $\sigma_3(G(3,7))\simeq G(3,7)^*$\\
\hline
$e_{123}+e_{456}+e_{147}+e_{257}+e_{367}+e_{367}$ & $34$ & $\PP^{34}$\\
\hline
\end{tabular}
\end{table}

\begin{rem}
 The variety $\sigma_{2,+}(G(3,n))$ is defined as the set of chords $\PP_{xy}^1$ such that $x,y\in G(3,n)$ and the corresponding 3-planes 
 intersect along a line. In the three tables we have $\text{Sing}(X^*)\simeq \sigma_{2,+}(G(3,n))^*$. 
In Table \ref{t6} this variety is self-dual.  The variety $\sigma_{3,+}(G(3,8))$, in Table \ref{t8} is defined as the closure of 
the set of planes $\PP_{xyz}^2$ passing through $x,y,z\in G(3,8)$ such that the three corresponding $3$-planes in $\Lambda^3 \CC^8$ meet along a line.
\end{rem}

\begin{table}[h]
\caption{$SL_8$-orbits in $\PP(\Lambda^3\CC^8)$}\label{t8}
\begin{tabular}{|c|c|c|c|}
  \hline
\text{Orbit} & \text{Representative} &\text{Dimension} & \text{Geometric interpretation} \\
\hline
II &$e_{123}$ & $15$ & $G(3,8)$\\
\hline
III & $e_{123}+e_{145}$ & $24$ & $\sigma_{2,+}(G(3,8))$\\
\hline
IV &$e_{124}+e_{135}+e_{236}$ &$30$ & $\tau(G(3,8))$\\
\hline
V &$e_{123}+e_{456}$ & $31$ & $\sigma_2(G(3,8))$\\
\hline
VI & $e_{123}+e_{145}+e_{167}$ & $27$ & $\sigma_{3,+}(G(3,8))$ \\
\hline
VII& $e_{125}+e_{136}+e_{147}+e_{234}$ & $34$& $X_7$\\
\hline
VIII& $e_{134}+e_{256}+e_{127}$ &$37$ &  $J(G(3,8),\tau(G(3,8)))^*$\\
\hline
IX& $e_{125}+e_{346}+e_{137}+e_{247}$ & $40$&  $J(G(3,8),\sigma_{+}(G(3,8)))^*$\\
\hline
X & $e_{123}+e_{456}+e_{147}+e_{257}+e_{367}$ & $41$ & $\sigma_{3,+}(G(3,8))^*$\\
\hline
XI & $e_{127}+e_{138}+e_{146}+e_{235}$ & $39$ & $T(G(3,8),\sigma_{2,+}(G(3,8)))$\\
\hline
XII & $e_{128}+e_{137}+e_{146}+e_{236}+e_{245}$ & $42$ & $T(G(3,8),\tau(G(3,8)))^*$ \\
\hline
XIII & $e_{135}+e_{246}+e_{147}+e_{238}$& $43$ & $X_{13}\simeq X_{13}^*$\\
\hline
XIV & $e_{138}+e_{147}+e_{156}+e_{235}+e_{246}$ &$45$  &$T(G(3,8),\tau(G(3,8)))$ \\
\hline
XV & $e_{128}+e_{137}+e_{146}+e_{247}+e_{256}+e_{345}$ & $47$ &  $T(G(3,8),\sigma_{2,+}(G(3,8)))^*$\\
\hline
XVI & $e_{156}+e_{178}+e_{234}$ & $40$ & $J(G(3,8),\sigma_{+}(G(3,8)))$\\
\hline
XVII & $e_{158}+e_{167}+e_{234}+e_{256}$ & $46$  & $J(G(3,8),\tau(G(3,8)))$ \\
\hline
XVIII & $e_{148}+e_{157}+e_{236}+e_{245}+e_{347}$ & $49$ & $X_7 ^*$\\
\hline
XIX & $e_{134}+e_{234}+e_{156}+e_{278}$ &$47$  & $\sigma_2(G(3,8))^* \simeq \sigma_3(G(3,8))$\\
\hline
XX& $e_{137}+e_{237}+e_{256}+e_{148}+e_{345}$ &  $51$ & $\tau(G(3,8))^*$  \\
\hline
XXI & $e_{138}+e_{147}+e_{245}+e_{267}+e_{356}$ & $52$ & $\text{Sing}(G(3,8)^*)\simeq \sigma_{2,+}(G(3,8))^*$\\
\hline
XXII &$e_{128}+e_{147}+e_{236}+e_{257}+e_{358}+e_{456}$  & $54$ & $G(3,8)^*$  \\
\hline
XXIII& $e_{124}+e_{134}+e_{256}+e_{378}+e_{157}+e_{468}$& $55$ & $\PP^{55}$\\
\hline
\end{tabular}
\end{table}

\begin{rem}
 The varieties corresponding to orbits $VII$ and $XIII$ (notations of D. \v{Z}. Djokovi\'c \cite{Djo}) are denoted by $X_7$ and $X_{13}$ because we do not have a geometric interpretation for those orbits. 
The variety $X_{13}$ is a new example of non-smooth self-dual variety. For more examples of non-smooth self-dual varieties arising 
from group actions see \cite{Po,PoTe}.
\end{rem}

\proof The first five orbits are clearly identified from their representatives. The same is true for $\overline{\mathcal{O}}_{\text{XVI}}$ (closure of orbit XVI): 
its representative $e_{156}+e_{178}+e_{234}$ is a general point of $J(G(3,8),\sigma_{2,+}(G(3,8)))$. That variety has the expected dimension,\[dim(J(G(3,8),\sigma_{2,+}(G(3,8)))=dim(G(3,8))+dim(\sigma_{2,+}(G(3,8)))+1=40\] therefore there exists an orbit of dimension $39$ which corresponds 
to $T(G(3,8),\sigma_{2,+}(G(3,8))$ (this is a consequence of the Fulton-Hansen connectedness Theorem, \cite{Z2}). 
But there is only one orbit of dimension $39$, thus $\overline{\mathcal{O}}_{\text{XI}}= T(G(3,8),\sigma_{2,+}(G(3,8))$. 
The variety $J(G(3,8),\sigma_{2,+}(G(3,8))$ is included in $\sigma_3(G(3,8))$ and we know by \cite{AOP} that
$\sigma_3(G(3,8))$ has dimension $47$. The order among the orbits (\cite{Djo} and the appendix of this article) proves that $\overline{\mathcal{O}}_{\text{XIX}}=\sigma_3(G(3,8))$.
The representative of $\mathcal{O}_{\text{XVII}}$, $\underbrace{e_{158}+e_{167}+e_{256}}_{\in \tau(G(3,8))}+\underbrace{e_{234}}_{\in G(3,8)}$
 belongs to $J(G(3,8),\tau(G(3,8))$. Thus 
$\overline{\mathcal{O}}_{\text{XVII}}\subset J(G(3,8),\tau(G(3,8))$. But $\text{dim}(\mathcal{O}_{\text{XVII}})=46$ which is the expected
 dimension of $J(G(3,8),\tau(G(3,8))$. It proves $\overline{\mathcal{O}}_{\text{XVII}}=J(G(3,8),\tau(G(3,8))$.
Then there exists an orbit of dimension $45$ which corresponds by the Fulton-Hansen Theorem to $T(G(3,8),\tau(G(3,8))$ and this orbit is $\mathcal{O}_{\text{XIV}}$.

To prove the duality between the orbits we first identify $\Lambda^3\CC^8$ and $\Lambda^3(\CC^*)^8$ by the usual pairing 
$<e_{ijk}, e^{rst}>=det(e^u(e_l))_{l=i,j,k;u=r,s,t}$
where $e^1,\dots,e^8$ is a basis of $(\CC^*)^8$. For each orbit $\mathcal{O}$ we construct 
$y\in (T_{x} \mathcal{O})^\perp$ such that $y$ is a representative of an orbit $\mathcal{O}'$.
That construction shows $\overline{\mathcal{O}}'\subset \overline{\mathcal{O}}^*$. Then the order among the orbits allow to conclude:
\begin{enumerate}
 \item Clearly $G(3,8)^*\simeq \overline{\mathcal{O}}_{\text{XXII}}$.
\item Let $x=e_{846}+e_{857}$ a representative of $\sigma_{2,+}(G(3,8))$, the element $y=e_{138}+e_{147}+e_{245}+e_{267}+e_{356}\in (T_x \sigma_{2,+}(G(3,8)))^\perp$. 
Thus $\PP(\overline{G.y})\subset \sigma_{2,+}(G(3,8))^*$. But $y$ is a representative of $\mathcal{O}_\text{XXI}$.
Then we conclude $\overline{\mathcal{O}}_{\text{XXI}}\subset \sigma_{2,+}(G(3,8))^*\subsetneq \overline{\mathcal{O}}_{\text{XXII}}$ and therefore 
$\overline{\mathcal{O}}_\text{XXI}\simeq \sigma_{2,+}(G(3,8))^*$ because there is no orbit between $\mathcal{O}_{\text{XXI}}$ and $\mathcal{O}_{\text{XXII}}$.
\item In the next table we give for each orbit $\mathcal{O}$ a representative $x$, an element $y\in (T_{x}\overline{\mathcal{O}})^\perp$
and the orbit corresponding to $\PP(\overline{G.y})$. 
This table proves  $\PP(\overline{G.y})\subset\overline{\mathcal{O}}^*$ and we conclude to the equality by looking at the order among the orbits (see the appendix and \cite{Djo}):
\begin{table}[h]
\caption{Duality between the orbits}
 \begin{tabular}{|c|c|c|c|}
\hline
Orbit $\overline{\mathcal{O}}=\PP(\overline{G.x})$  & Representative $x$ & Representative $y\in (T_{x}\mathcal{O})^\perp$ & Orbit $\PP(\overline{G.y})$\\
\hline
$\sigma_{2,+}(G(3,8))$ & $e_{846}+e_{857}$ & $e_{138}+e_{147}+e_{245}+e_{267}$& $\mathcal{O}_{\text{XXI}}$\\
                      &                     & $+e_{356}$ 	&\\
\hline
$\tau(G(3,8)$ & $e_{467}+e_{368}+e_{578}$ &$e_{137}+e_{237}+e_{256}+e_{148}$ &$\mathcal{O}_{\text{XX}}$ \\
              &                           & $+e_{345}$ & \\
\hline 
$\sigma(G(3,8))$ & $e_{357}+e_{468}$ & $e_{134}+e_{234}+e_{156}+e_{278}$ & $\mathcal{O}_{\text{XIX}}$\\
\hline
$X_7$ & $e_{835}+e_{872}+e_{864}+e_{567}$ & $e_{148}+e_{157}+e_{236}+e_{245}$ & $\mathcal{O}_{\text{XVIII}}$\\
      &                                 &       $+e_{347}$ &\\
\hline
$\mathcal{O}_{\text{VIII}}$ & $e_{134}+e_{256}+e_{127}$ & $e_{832}+e_{851}+e_{764}+e_{735}$ & $J(G(3,8),\tau(G(3,8)))$\\
\hline
$\mathcal{O}_{\text{IX}}$ & $e_{125}+e_{346}+e_{137}+e_{247}$ & $e_{841}+e_{823}+e_{567} $ & $J(G(3,6),\sigma_+(G(3,8)))$\\
\hline
$T(G(3,8),\sigma_{2,+}(G(3,8)))$ & $e_{821}+e_{836}+e_{875}+e_{472}$ & $e_{128}+e_{137}+e_{146}+e_{247}$ & $\mathcal{O}_{\text{XV}}$\\
                                 &                                   &$+e_{256}+e_{345}$& \\
\hline
$\mathcal{O}_{\text{XII}}$ & $e_{128}+e_{137}+e_{146}+e_{236}$ & $e_{812}+e_{865}+e_{834}+e_{731}$ & $T(G(3,8),\tau(G(3,8))$\\
                          &       $+e_{245}$                    &        $+e_{754}$                        &\\
\hline
$X_{13}$ & $e_{752}+e_{861}+e_{763}+e_{845}$ & $e_{135}+e_{246}+e_{147}+e_{238}$ & $X_{13}$ \\
\hline
$\sigma_{3,+}(G(3,8))$ & $e_{815}+e_{826}+e_{834}$                  & $e_{123}+e_{456}+e_{147}+e_{257}$ & $\mathcal{O}_{X}$ \\
                      &                                               &       $+e_{367}$  &\\         
\hline
 \end{tabular}

\end{table}

\end{enumerate}

\begin{rem}
The representatives for orbits with geometric interpretation can easly be identified. For instance it is clear that $e_{467}+e_{368}+e_{578}$ is a representative of the tangential variety. 
The varieties for which we need to explicitly 
write the 
new indexation of the basis to identify the representative are $T(G(3,8),\sigma_{2,+}(G(3,8)))$, $X_{13}$, $T(G(3,8),\tau(G(3,8))$ and $X_7$.
For instance for $T(G(3,8),\tau(G(3,8))$ we consider the following change of basis $g$: $e_8\rightarrow e_1$, $e_1\rightarrow e_3$, $e_2\rightarrow e_8$, $e_5\rightarrow-e_4$, 
$e_6\rightarrow e_7$, $e_3\rightarrow e_5$, $e_4\rightarrow e_6$, $e_7\rightarrow -e_2$, with
 $g\in SL_8$ and $g(e_{812}+e_{865}+e_{834}+e_{731}+e_{754})=e_{138}+e_{147}+e_{156}+e_{235}+e_{246}$  which is the representative of $T(G(3,8),\tau(G(3,8))$ in Table \ref{t8}. Similar changes of basis exist for the remaining cases.
\end{rem}

\setcounter{theo}{1}
\section{Veronese embeddings and Segre products}\label{hyper}
The Proposition \ref{crit} can be  used to get similar results on Veronese re-embeddings and Segre products of smooth projective varieties.
\begin{theo}\label{veroseg}
 Let $X\subset \PP^n$ and $Y\subset \PP^m$ two smooth nondegenerate projective varieties. Then 
\begin{enumerate}
 \item For $d\geq 2$, $\sigma_2(v_d(X))^*$ is a codimension one subvariety of $v_d(X)^*$ if and only if $(X,d)\neq (\PP^n,2)$.
\item $\sigma_2(X\times Y)^*$ is a codimension one subvariety of $(X\times Y)^*$ when either 
$X\times Y\neq X\times \PP^m$ or $\sigma_2(X)$ is not defective. 
\end{enumerate}
\end{theo}

\proof We calculate the tangent space and the second osculating space:
\begin{enumerate}
 \item Let $(x^d,y^d)$ be a general pair of points of $v_d(X)$ then using Leibniz's rule we have $\hat{T}^{(2)}_{x^d} X=\hat{T}^{(2)} _x X\circ x^{d-1}+\hat{T}_x X\circ \hat{T}_x X\circ x^{d-2}$ and 
$\hat{T}_{y^d} v_d(X)=\hat{T}_y X\circ y^{d-1}$. The intersection $\hat{T}^{(2)}_{x^d} X\cap \hat{T}_{y^d} v_d(X)\neq \{0\}$ if and only if $d=2$ and 
$y\in \hat{T}_x X$ i.e. $X=\PP^n$. Thus Proposition \ref{crit} applies.
On the other hand it is known that $v_2(\PP^n)^*$ is regular in codimension one (\cite{Z2}).
\item Let $(x\otimes y,u\otimes v)$ be a general pair of points of $X\times Y$. Then $\hat{T}^{(2)}_{x\otimes y} (X\times Y)=\hat{T}^{(2)} _x X\otimes y+\hat{T}_x X\otimes \hat{T}_y Y+x\otimes \hat{T}^{(2)} _y Y$
and $\hat{T}_{u\otimes v} (X\times Y)=\hat{T}_u X\otimes v+u\otimes \hat{T}_v Y$. 
The intersection $\hat{T}^{(2)}_{x\otimes y} (X\times Y)\cap \hat{T}_{u\times v} (X\times Y)$ does not reduced to $0$ only if 
\begin{enumerate}
 \item $\tilde{T}_x X=\PP^n$ and $\hat{T}_y Y\cap \hat{T}_w Y\neq \{0\}$ i.e. $X=\PP^n$ and $\sigma_2(Y)$ is defective.
 \item $\tilde{T}_u X=\PP^n$ and $\tilde{T}^{(2)} _y Y =\PP^m$ i.e. $X=\PP^n$ and $\sigma_2(Y)$ is defective.
\item $\tilde{T}_w Y=\PP^m$ and $\tilde{T}^{(2)} _x X=\PP^n$ i.e. $Y=\PP^m$ and $\sigma_2(X)$ is defective.
\item $\tilde{T}_y Y=\PP^m$ and $\hat{T}_x X \cap \hat{T}_u X\neq \{0\}$ i.e. $Y=\PP^m$ and $\sigma_2(X)$ is defective.
\end{enumerate}
Thus Proposition \ref{crit} applies outside the previous four cases. $\Box$
\end{enumerate}

\par{\bf Back to Hyperdeterminants:}
The steps we followed in $\S\ref{grass}$ allow us to recover the result on normality of hyperdeterminants. Let $X=\PP^{k_1}\times \PP^{k_2}\times\dots \PP^{k_s}\subset \PP^{(k_1+1)(k_2+1)\dots(k_s+1)-1}$. Suppose $k_i\leq k_1+k_2+\dots+\hat{k_i}+\dots+k_s$ so that $X^*$ is a hypersurface (\cite{GKZ1}) :
\begin{enumerate}
\item The second part of Theorem \ref{veroseg}  shows that for hyperdeterminants the only chance to get a dual variety regular in codimension one is when
we consider $X=\PP^{k_1}\times \PP^{k_2}\times \PP^{k_3}$. 
\item Similar arguments to $\S\ref{grass}$ prove that  $(\PP^{k_1}\times \PP^{k_2}\times \PP^{k_3})^*$ is singular in codimension 
$1$ when $k_1+k_2+k_3\geq 6$. More precisely the calculation on the rank of specific quadrics of 
$|II_{\PP^{k_1}\times \PP^{k_2}\times \PP^{k_3}}|$ leads to consider the matrices of type 
$\begin{pmatrix}
     0 & A & B \\
^t A & 0 & C\\
^tB & ^t C & 0
 \end{pmatrix}$ with blocks $A,B,C$ respectively of size $k_1\times k_2$, $k_1\times k_3$, and $k_2\times k_3$ and 
with the additional condition $a_{11}=b_{11}=c_{11}=0$. The condition on the corner entry of each block appears from the same reason as
 condition $(\star\star)$ in section $\S\ref{grass}$. 
\item To finish the proof we consider the following orbits:
\begin{enumerate}
 \item action of $SL_3\times SL_3\times SL_2$ on $\CC^3\otimes \CC^3\otimes \CC^2$, 
\item action of $SL_3\times SL_2\times SL_2$ on $\CC^3\otimes \CC^2\otimes \CC^2$,
\item action of $SL_2\times SL_2\times SL_2$ on $\CC^2\otimes \CC^2\otimes \CC^2$.
\end{enumerate}
All of those group actions have finitely many orbits and there is bijection between the orbits in $\PP(V_1\otimes V_2\otimes V_3)$ and $\PP(V_1 ^*\otimes V_2^*\otimes V_3 ^*)$ (see \cite{Par}). 
It follows that we find a hypersurface regular in codimension 1 only in (c).
\end{enumerate}

\appendix
\section{Orbits decompositions}
We reproduce the order relation among the orbits for $k=3$ and $n=6,7,8$. 
We use the graphical notation coming from \cite{Shou} and \cite{AOP}. Each node corresponds to an element of
the basis of $\CC^n$. The linked nodes represent a tri-vector and for each diagram a representative of an orbit is the sum of the corresponding trivectors. We number the nodes only for the orbit corresponding to the ambient space. 
The representatives are those of Gurevich's book \cite{Gu}.
\begin{figure}[!h]
\begin{minipage}{6cm}
\begin{center}
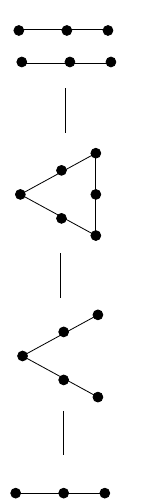
\end{center}
\caption{Orbits in $\PP(\Lambda^3\CC^6)$}
\end{minipage}
\begin{minipage}{6cm}
\begin{center}
 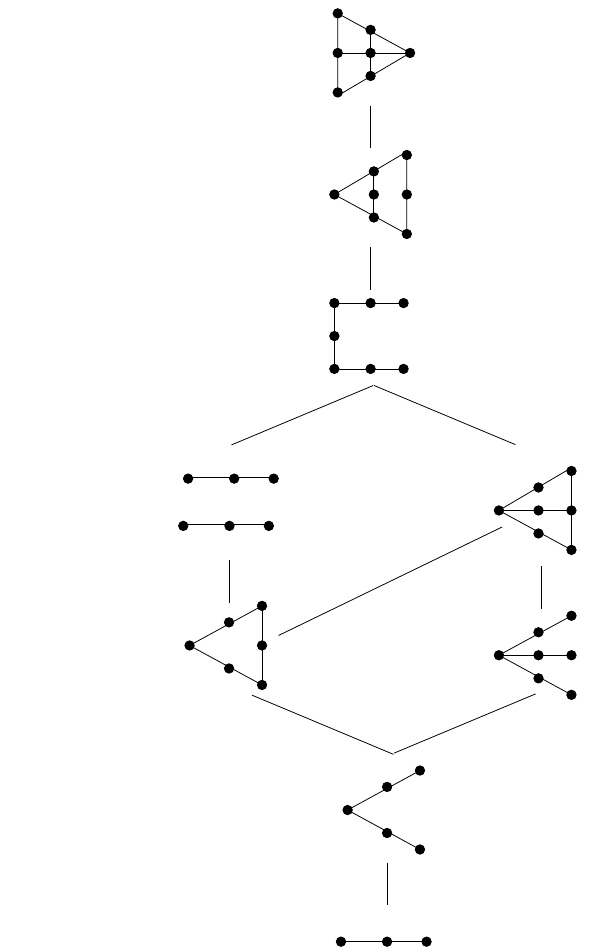
\end{center}
\caption{Orbits in $\PP(\Lambda^3\CC^7)$}
\end{minipage}
\end{figure}

\begin{figure}[!h]
\begin{center}
 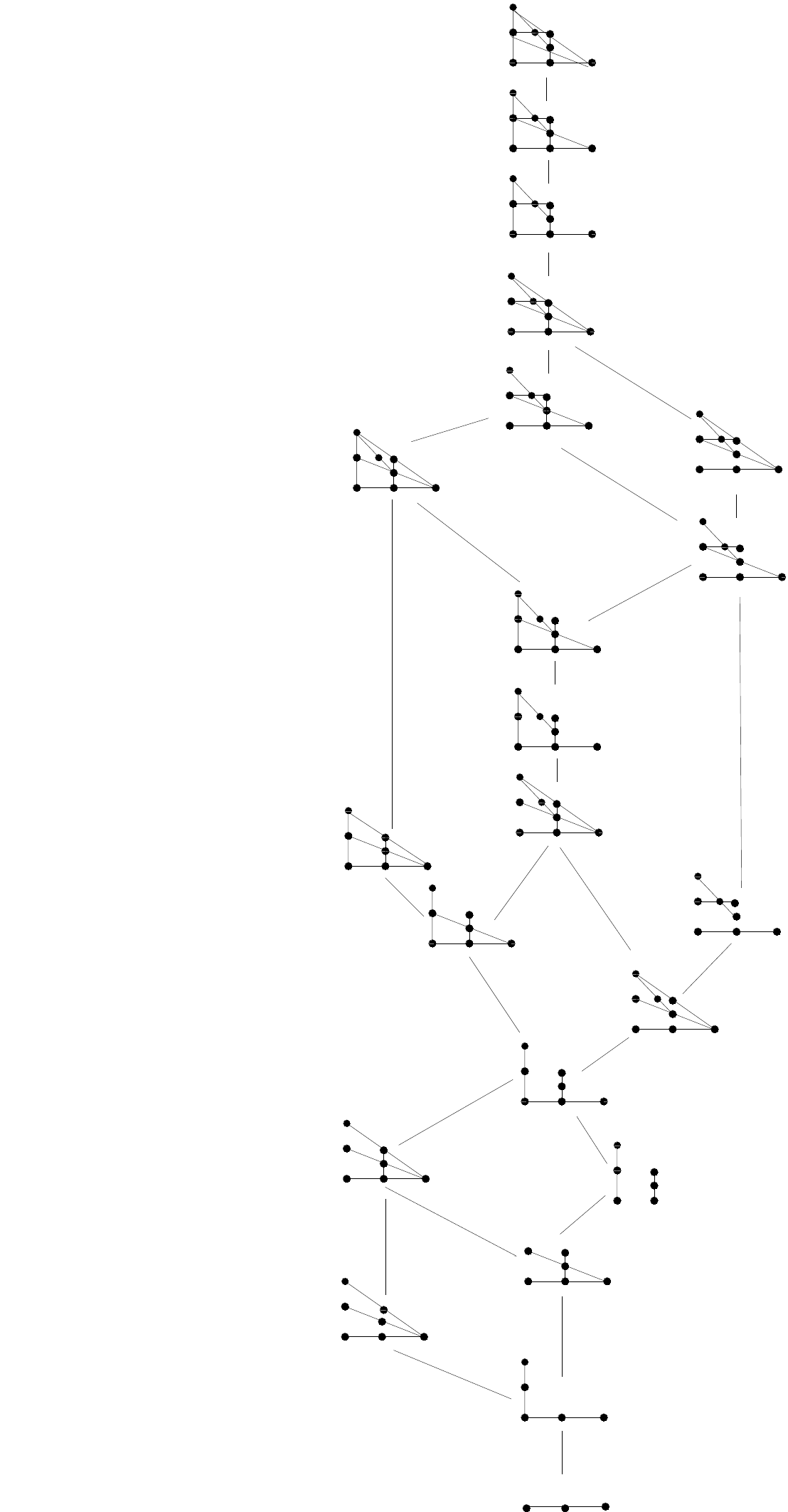
\end{center}
\caption{Orbits in $\PP(\Lambda^3 \CC^8)$}
\end{figure}
\end{document}